\documentclass[11pt]{amsart}
\usepackage{amsmath,amsthm,amssymb, amscd, amsfonts}
\usepackage[all]{xy}
\usepackage{leftidx}
\usepackage[latin1]{inputenc}
\usepackage{enumerate}
\usepackage[dvipsnames]{xcolor}

\theoremstyle{plain}

\newtheorem{theorem}{Theorem}[section]
\newtheorem{corollary}[theorem]{Corollary}

\newtheorem{proposition}[theorem]{Proposition}
\newtheorem{lemma}[theorem]{Lemma}

\theoremstyle{definition}

\theoremstyle{remark}
\newtheorem{remark}[theorem]{Remark}

\numberwithin{equation}{section}\theoremstyle{plain}

\renewcommand{\1}{\textbf{1}}

\newcommand{\A}{{\mathcal A}}
\newcommand{\B}{{\mathcal B}}
\newcommand{\C}{{\mathcal C}}
\newcommand{\D}{{\mathcal D}}

\newcommand{\Z}{{\mathcal Z}}

\newcommand{\E}{{\mathcal E}}
\newcommand{\U}{{\mathcal U}}

\newcommand{\Rep}{\operatorname{Rep}}

\newcommand\Aut{\operatorname{Aut}}
\newcommand\Irr{\operatorname{Irr}}
\newcommand\FPdim{\operatorname{FPdim}}

\newcommand\vect{\operatorname{Vec}}
\newcommand\svect{\operatorname{sVec}}
\newcommand\id{\operatorname{id}}

\newcommand\Tr{\operatorname{Tr}}

\newcommand\Hom{\operatorname{Hom}}
\newcommand\Fun{\operatorname{Fun}}

\begin{document}
\title{Classification of certain weakly integral fusion categories}
\author{Jingcheng Dong}

\address[Dong]{College of Mathematics and Statistics, Nanjing University of Information Science and Technology, Nanjing 210044, China}
\email{jcdong@nuist.edu.cn}



\date{\today}

\begin{abstract}
We prove that braided fusion categories of Frobenius-Perron $p^mq^nd$ or $p^2q^2r^2$ are weakly group-theoretical, where $p,q,r$ are distinct prime numbers, $d$ is a square-free natural number such that $(pq,d)=1$. As an application, we obtain that weakly integral braided fusion categories of Frobenius-Perron dimension less than $1800$ are weakly group-theoretical, and weakly integral braided fusion categories of odd dimension less than $33075$ are solvable. For the general case, we prove that fusion categories  (not necessarily braided) of Frobenius-Perron dimension $84$ and $90$ either solvable or group-theoretical. Together with the results in the literature, this shows that every weakly integral fusion category of Frobenius-Perron dimension less than $120$ is either solvable or group-theoretical. Thus we complete the classification of all these fusion categories in terms of Morita equivalence.
\end{abstract}

\maketitle

\section{Introduction}\label{sec1}
 A fusion category $\C$ is a $k$-linear semisimple rigid tensor category with finitely many isomorphism classes of simple objects, finite-dimensional vector spaces of morphisms and such that the unit object $\1$ is simple. The theory of fusion category arises from many areas of mathematics and physics, such as semisimple Hopf algebras, quantum groups \cite{BaKi2001lecture}, vertex operator algebras \cite{2016DongWang} and topological quantum field theory \cite{Turaer1994}.

The notion of a weakly group-theoretical fusion category is introduced in \cite{etingof2011weakly}. By definition, a fusion category $\C$ is called weakly group-theoretical if it is Morita equivalent to a nilpotent fusion category. In particular, every weakly group-theoretical fusion category is weakly integral, i.e., it has integral Frobenius-Perron dimension. This fact motivates the conjecture that every weakly integral fusion category is weakly group-theoretical, see \cite{etingof2011weakly}. Some examples of weakly group-theoretical fusion categories were obtained in \cite{etingof2011weakly,Dong2019,Dongpointed2020,DongChenWang2022,DNS2019,DONG2021386}. Natale also obtained some examples under the assumption of non-degeneracy, see \cite{natale2013weakly}. In fact, all known  weakly integral fusion categories are weakly group-theoretical.

\medbreak
People's interest in weakly group-theoretical fusion categories also comes from the fact that they have strong Frobenius property \cite[Theorem 1.5]{etingof2011weakly}. That is, for every weakly group-theoretical fusion category $\C$ and every indecomposable $\C$-module category $\mathcal{M}$ and any simple object $X$ in $\mathcal{M}$ the number   $\frac{\FPdim(\C)}{\FPdim(X)}$ is an algebraic integer, where $\FPdim$ stands for the Frobenius-Perron dimension. If we take $\mathcal{M}=\C$ then the strong Frobenius property of a fusion category implies the usual Frobenius property which was  conjectured by  Kaplansky in the setting of semisimple Hopf algebras \cite{Kaplansky1975}.

\medbreak
The class of solvable or group-theoretical fusion categories are special cases of weakly group-theoretical fusion categories.  They are both interesting examples which come from finite group theory. Moreover, the later one can be completely classified by  finite groups and their cohomology \cite{etingof2005fusion}. So it is an interesting task to determine which class of weakly group-theoretical fusion categories are solvable or group-theoretical.

\medbreak
The present paper is devoted to extending the results obtained in \cite{etingof2011weakly} and discarding the assumption of non-degeneracy in \cite{natale2013weakly}. Our first result is listed below.

\begin{theorem}
Let $\C$ be a braided fusion category of Frobenius-Perron dimension $p^mq^nd$ or $p^2q^2r^2$, where $p,q,r$ are distinct prime numbers, $d$ is a square-free natural number such that $(pq,d)=1$. Then $\C$ is weakly group-theoretical.
\end{theorem}

Subsequently, all weakly integral braided fusion categories of dimension less than $1800$ are weakly group-theoretical. In particular, weakly integral braided fusion categories of odd dimension less than $33075$ are solvable.

\medbreak
Our second result is the following theorem.
\begin{theorem}
Fusion categories (not necessarily braided) of Frobenius-Perron dimension $84$ and $90$ are solvable or group-theoretical.
\end{theorem}

Our conclusion, together with results in the literature, shows that every weakly integral fusion category of Frobenius-Perron dimension less than $120$ is either solvable or group-theoretical.

\medbreak
The paper is organized as follows. In Section \ref{sec2}, we recall some basic definitions and results which will be used throughout. In Section \ref{sec3}, we study braided fusion category of Frobenius-Perron $p^mq^nd$ and $p^2q^2r^2$.
In Section \ref{sec41}, we study the existence of nontrivial symmetric categories in the Drinfeld center of a fusion category. Fusion categories of dimension $84$ and $90$ will be studied in Section \ref{sec42} and \ref{sec43}, respectively.

\medbreak
Throughout this paper, we shall work over an algebraically closed field $k$ of characteristic $0$. We refer to \cite{egno2015} for the main notions about fusion categories.

\section{Preliminaries}\label{sec2}
\subsection{Frobenius-Perron dimensions}\label{sec2.1}
Let $\C$ be a fusion category and let $\Irr(\C)$ denote the set of isomorphism classes of simple objects in $\C$. Then $\Irr(\C)$ is a basis of the Grothendieck ring $K_0(\C)$ of $\C$. The Frobenius-Perron dimension $\FPdim(X)$ of $X \in \Irr(\C)$ is defined as the largest eigenvalue of
the matrix of left multiplication by $X$ in the Grothendieck ring with respect to the basis. The Frobenius-Perron dimension of $\C$ is the number
$$\FPdim (\C)=\sum_{X \in \Irr(\C)} (\FPdim X)^2.$$

\medbreak
A fusion category $\C$ is called weakly integral if $\FPdim(\C)$ is an integer. A fusion category $\C$ is  called integral if $\FPdim(X)$  is an integer for every $X\in\Irr(\C)$.

\medbreak
A simple object $X\in \C$ is called invertible if $\FPdim(X)=1$. A pointed fusion category is a fusion category $\C$ whose  simple objects are all invertible.  If $\C$ is a pointed fusion category, then $\C$ is equivalent to the category of $G$-graded vector spaces with associativity constraint given by a $3$-cocycle $\omega\in H^3(G,k^{\times})$, denoted by  $\vect_{G,\omega}$. This fact is originally due to Eilenberg and Mac Lane. See also \cite{frohlich1993} for details.

We use $\C_{pt}$ to denote the largest pointed fusion subcategory in $\C$. Let $G(\C)$ be the group of isomorphism classes of invertible simple objects of a fusion category $\C$. Then $|G(\C)|=\FPdim(\C_{pt})$.

\medbreak
Let $Y $ be an object of $\C$ and write $Y = \sum_{X \in \Irr
(\C)} m(X, Y) X$, where $m(X, Y) \in \mathbb Z$. The  integer
$m(X, Y)$ is called the multiplicity of $X$ in $Y$. By \cite[Theorem
9,10]{Nichols1996}, we have the following results which will be used in our computation.

Let $\C$  be an integral fusion category and let $X, Y, Z$ be objects of $\C$. Then we have
$m(X,Y)=m(X^*,Y^*)$, and
$$m(X,Y\otimes Z)=m(Y^*,Z\otimes X^*)=m(Y,X\otimes Z^*).$$

Let $g$ be an element in $G(\C)$.  Then we have
$m(g,X\otimes Y)=1$ if and only if $Y=X^*\otimes g$, otherwise it is $0$. In
particular, $m(g,X\otimes Y)=0$ if $\FPdim X\neq \FPdim Y$. Let $X\in
\Irr(\C)$. Then for all $g \in G(\C)$, $m(g,X\otimes X^{*})>0$ if and
only if $m(g,X\otimes X^{*})= 1$ if and only if $g\otimes X=X$. The set of
isomorphism classes of such invertible objects will be denoted by
$G[X]$. Thus $G[X]$ is a subgroup of $G(\C)$ of order dividing
$(\FPdim X)^2$. In particular, for all $X\in \Irr(\C)$, we have a decomposition
\begin{equation}\label{second}
\begin{split}
X\otimes X^*=\bigoplus_{g\in G[X]}g\oplus\sum_{Y\in \Irr(\C)-G[X]}  m(Y,X\otimes X^*)Y.
\end{split}
\end{equation}

It is known that the group $G(\C)$ acts on the set $\Irr(\C)$ by left tensor
multiplication. This action preserves Frobenius-Perron dimensions
and, for $X\in \Irr(\C)$, $G[X]$ is the stabilizer of $X$ in
$G(\C)$.

\medbreak
In fact, these results  and Theorem \ref{thm-dim-2} below were established in the case
where $\C$ is the category of finite-dimensional representations
of a semisimple Hopf algebra. Because their proofs only make use of
the properties of the Grothendieck ring, these proofs also work  \textit{mutatis
mutandis} in the fusion category setting, and thus we
omit its proof.

\medbreak
Let $1=d_0, d_1,\cdots, d_s$, $s \geq 0$, be positive
real numbers such that $1 = d_0 < d_1< \cdots < d_s$, and let
$n_1,n_2,\cdots,n_s$ be positive integers. We shall
say that $\C$ is \emph{of type} $(d_0,n_0;
d_1,n_1;\cdots;d_s,n_s)$ if, for all $i = 0, \cdots, s$, $n_i$ is
the number of the non-isomorphic simple objects of
Frobenius-Perron dimension $d_i$.

The following theorem is a restatement of \cite[Theorem
11]{Nichols1996} in the context of fusion categories.

\begin{theorem}\label{thm-dim-2}
Let $\C$ be an integral fusion category and let $X$ be a simple object of dimension $2$. Then at least one of the following holds:

(1)\, $G[X] \neq \1$.

(2)\, $\C$ has a fusion subcategory $\D$ of type $(1, 2; 2, 1; 3, 2)$, such that $X \notin \Irr(\D)$, and $\D$ has an invertible object $g$
of order $2$ such that $g\otimes X\neq X$.

(3)\, $\C$ has a fusion subcategory of type $(1, 3; 3, 1)$ or $(1, 1; 3, 2; 4, 1; 5, 1)$.

In particular, if $G[X] = \1$, then $\C$ contains a fusion subcategory of dimension  $12$, $24$, or $60$.
\end{theorem}

\subsection{Extensions of a fusion category}\label{sec2.2}
Let $G$ be a finite group. A $G$-grading of $\C$ is a decomposition of $\C$ as a direct sum of full Abelian subcategories $\C=\oplus_{g\in G}\C_g$, such that $(\C_g)^*=\C_{g^{-1}}$ and the tensor product $\otimes:\C\otimes\C\to\C$ maps $\C_g\times \C_h$ to $\C_{gh}$.  The grading $\C=\oplus_{g\in G}\C_g$ is called faithful if $\C_g\neq0$ for all $g\in G$. We say that $\C$ is an $G$-extension of $\D$ if the grading is faithful and the trivial component is $\D$.

By \cite[Theorem 3.10]{gelaki2008nilpotent}, every weakly  integral fusion category is a $G$-extension of an integral fusion category $\D$, where $G$ is an elementary abelian $2$-group.

Let $\C$ be a fusion category. The fusion subcategory of $\C$ generated by simple objects in $X\otimes X^{*}$ for all $X\in\Irr(\C)$ is called the adjoint subcategory of $\C$ and is denoted by $\C_{ad}$. By \cite[Corollary 3.7]{gelaki2008nilpotent} every fusion category has a canonical faithful grading $\C=\oplus_{g\in \U(\C)}\C_g$ with trivial component $\C_{ad}$. This grading is called the universal grading of $\C$ and $\U(\C)$ is called the universal grading group of $\C$.

Let $\C=\oplus_{g\in G}\C_g$ be a $G$-extension of $\D$. Then $\FPdim(\C_g)=\FPdim(\C_h)$ for all $g,h\in G$ and $\FPdim(\C)=|G|\FPdim(\D)$, see \cite[Proposition 8.20]{etingof2005fusion}.

\subsection{Weakly group-theoretical fusion categories}
Let $\C$ be a fusion category
and let $\mathcal{M}$ be an indecomposable right $\C$-module category. Let $C^*_{\mathcal{M}}$ denote the category of $\C$-module endofunctors of $\mathcal{M}$. Then  $C^*_{\mathcal{M}}$ is a fusion category, called the dual of $\C$ with respect
to $\mathcal{M}$ \cite{etingof2005fusion,ostrik2003module}. Two fusion categories $\C$ and $\D$ are Morita equivalent if $D$ is equivalent to $C^*_{\mathcal{M}}$ for some indecomposable right $\C$-module category $\mathcal{M}$.

A fusion category $\mathcal{\C}$ is said to be (cyclically) nilpotent if there is a sequence of
fusion categories
\begin{equation}\label{nilp}
\begin{split}
\mathcal{\C}_0 ={\rm Vec}, \mathcal{\C}_1, \cdots, \mathcal{\C}_n =\mathcal{\C}
\end{split}
\end{equation}
and a sequence of finite (cyclic) groups $G_1,\cdots, G_n$ such that $\mathcal{\C}_i$ is obtained from $\mathcal{\C}_{i-1}$ by a $G_i$-extension.

A fusion category is called weakly group-theoretical if it is Morita equivalent to a nilpotent fusion category. A fusion category is called group-theoretical if it is Morita equivalent to a pointed fusion category. A fusion category is called solvable if it is Morita equivalent to a cyclically nilpotent fusion category.

For weakly group-theoretical fusion categories, there is a version of the Feit-Thompson theorem: if
a weakly group-theoretical fusion category is odd-dimensional then it is solvable, see \cite[Proposition 7.1]{natale2012fusion}. The proposition below shows that similar result also holds for weakly group-theoretical fusion categories of dimension $2n$, where $n$ is odd.

Let us recall an old result in group theory before giving the proof of the proposition below. Let $G$ be a finite group. If $|G|$ is odd or $2n$ then $G$ is solvable,  where $n$ is odd, see  e.g., \cite[Theorem 1.35]{Isaacs}.

\begin{proposition}\label{wGt-sol}
Let $\C$ be a weakly group-theoretical fusion category. Assume that $\FPdim(\C)$ is $2n$, where $n$ is odd. Then $\C$ is solvable. 
\end{proposition}
\begin{proof}
By the definition above, $\C$ is Morita equivalent to a nilpotent fusion category $\D$. If $\D$ is pointed then $\D=\vect_{G,\omega}$ for a finite group $G$  with $|G|=\FPdim(\C)$ and a $3$-cocycle $\omega$. Since $G$ is solvable,  $\D$ is solvable and so is  $\C$,  by \cite[Proposition 4.5]{etingof2011weakly}. In fact, $\C$ is also group-theoretical by the definition of a group-theoretical fusion category.

\medbreak
We then assume that $\D$ is not pointed and hence there is a sequence of fusion categories
\begin{equation}\label{nilp2}
\begin{split}
\mathcal{\D}_0 ={\rm Vec}, \mathcal{\D}_1, \cdots, \mathcal{\D}_n =\mathcal{\D}
\end{split}
\end{equation}
and a sequence of finite groups $G_1,\cdots, G_n$ such that $\mathcal{\D}_i$ is obtained from $\mathcal{\D}_{i-1}$ by a $G_i$-extension.

It is clear that the order of $G_i$ divides $\FPdim(\C)$ and hence $G_i$ is solvable,  $1\leq i\leq n$. Moreover, the fusion subcategory $\D_1$ must be pointed and hence $\D_1=\vect_{K,\omega}$, where $K$ is a finite group with order dividing $\FPdim(\C)$. By \cite[Proposition 4.5(ii)]{etingof2011weakly}, $\D_1$ is solvable. Then $\D_i$ in the sequence is solvable by \cite[Proposition 4.5(ii)]{etingof2011weakly}. Again by \cite[Proposition 4.5(ii)]{etingof2011weakly}, $\C$ is solvable.
\end{proof}

\begin{corollary}\label{dimpqr}
Let $p<q<r$ be distinct prime numbers, and let $\C$ be a  fusion category of dimension $pqr$. Then $\C$ is solvable.
\end{corollary}
\begin{proof}
If $\C$ is integral then $\C$ is group-theoretical by \cite[Theorem 9.2]{etingof2011weakly}. Hence $\C$ is solvable by Theorem \ref{wGt-sol}.

If $\C$ is not integral then $\C=\oplus_{g\in\mathbb{Z}_2}\C_g$ by \cite[Theorem 3.10]{gelaki2008nilpotent}, where $\C_e$ is a fusion category of dimension $qr$. Hence $\C$ is solvable by \cite[Proposition 4.5]{etingof2011weakly}.
\end{proof}

\subsection{Braided fusion categories}\label{sec23}
A fusion category $\C$ is called braided if it admits a braiding $c$, where the braiding $c$ is a family of natural isomorphisms: $c_{X,Y}$:$X\otimes Y\rightarrow Y\otimes X$ satisfying the hexagon axioms for all $X,Y\in\C$ \cite{kassel1995quantum}.

\medbreak
Let $\C$ be a braided fusion category, and $\D\subseteq \C$ be a fusion subcategory. The M\"{u}ger centralizer $\D'$ of $\D$ in $\C$ is the category of all objects $Y\in \C$ such that $c_{Y,X}c_{X,Y}=\id_{X\otimes Y}$ for all $X\in \D$. The centralizer $\D'$ is again a fusion subcategory of $\C$. The M\"{u}ger center of $\C$ is the M\"{u}ger centralizer $\mathcal{Z}_2(\C):=\C'$ of $\C$ itself. A braided fusion category $\C$ is called nondegenerate if its M\"{u}ger center $\mathcal{Z}_2(\C)$ is trivial and it is called slightly degenerate if its M\"{u}ger center $\mathcal{Z}_2(\C)$ is equivalent to the category $\svect$ of super vector spaces.

\medbreak
Let $\C$ be a fusion category. The Drinfeld  center $\Z(\C)$ of $\C$  is defined as the category whose objects are pairs $(X, c_{-,X} )$, where $X$ is an object of $\C$ and $c_{-,X}$ is a natural family of isomorphisms $c_{V,X} : V\otimes X\to X\otimes V$, $V\in \C$, satisfying certain compatibility conditions, see \cite[Definition XIII.4.1]{kassel1995quantum}.  It is shown in \cite[Theorem 2.15, Proposition 8.12]{etingof2005fusion} that $\Z(\C)$ is a braided fusion category and $\FPdim(\Z(\C))=\FPdim(\C)^2$. In addition, $\Z(\C)$ is non-degenerate by \cite[Corollary 3.9]{drinfeld2010braided}.

\medbreak
A braided fusion category $\C$ is called symmetric if $c_{Y,X}c_{X,Y}=\id_{X\otimes Y}$ for all objects $X,Y\in\C$. A symmetric fusion category $\C$ is said to be Tannakian if it is equivalent to $\Rep(G)$ for some finite group $G$ as symmetric categories.

Let $G$ be a finite group and let $u\in G$ be a central element of order 2. Then the category $\Rep(G)$ has a braiding $c^u_{X,Y}$ as follows: for all\, $x\in X, y\in Y$,
\begin{equation}\label{eq2}
\begin{split}
c^u_{X,Y}(x\otimes y)=(-1)^{mn}y\otimes x\,\, \mbox{if}\,\,ux=(-1)^mx,uy=(-1)^ny.
\end{split}\nonumber
\end{equation}

Let $\Rep(G,u)$ be the fusion category $\Rep(G)$ equipped with the new braiding $c^u_{X,Y}$. Deligne proved that any symmetric fusion category is equivalent to some $\Rep(G,u)$ \cite{deligne1990categories}.

\begin{lemma}\cite[Corollary 2.50]{drinfeld2010braided}\label{lem22}
Let $\C$ be the symmetric fusion category $\Rep(G,u)$. Then one of the following holds:

(1) $\C$ is a Tannakian category;

(2) $\Rep(G/\langle u \rangle) \subseteq \C$ is a maximal Tannakian subcategory with dimension $\frac{1}{2}\FPdim(\C)$.

In particular, if $\FPdim(\C)$ is bigger than 2, then $\C$ always has a nontrivial Tannakian subcategory.
\end{lemma}


The following theorem is known as the M\"{u}ger Decomposition Theorem since it is due to M\"{u}ger \cite[Theorem 4.2]{muger2003structure} when $\C$ is modular.

\begin{theorem}{\cite[Theorem 3.13]{drinfeld2010braided}}\label{MugerThm}
Let $\C$ be a braided fusion category and let $\D$ be a non-degenerate subcategory of $\C$. Then $\C$ is braided equivalent to $\D\boxtimes \D'$, where $\D'$ is the centralizer of $\D$ in $\C$.
\end{theorem}

\subsection{Equivariantizations and braided $G$-crossed fusion categories}\label{subsec24}
Let $G$ be a finite group and $\mathcal{\C}$ be a fusion category. Let $\underline{G}$ denote the monoidal category whose objects are elements of $G$, morphisms are identities and tensor product is given by the multiplication in $G$. Let ${\rm\underline{Aut}}_{\otimes}\mathcal{\C}$ denote the monoidal category whose objects are tensor autoequivalences of $\C$, morphisms are isomorphisms of tensor functors and tensor product is given by the composition of functors.

\medbreak
An action of $G$ on $\mathcal{\C}$ is a monoidal functor
$$T:\underline{G}\to {\rm\underline{Aut}}_{\otimes}\mathcal{\C},\quad g\mapsto T_g$$
with the isomorphism $f^X_{g,h}: T_g(X)\otimes T_h(X)\cong T_{gh}(X)$, for every $X$ in $\mathcal{\C}$.

\medbreak
Let $\mathcal{\C}$ be a fusion category with an action of $G$. Then the fusion category $\mathcal{\C}^G$, called the $G$-equivariantization of $\mathcal{\C}$, is defined as follows \cite{bruguieres2000categories,drinfeld2010braided,muger2004galois}:

(1)\quad An object in $\mathcal{\C}^G$ is a pair $(X, (u^X_g)_{g\in G})$, where $X$ is an object of $\mathcal{\C}$ and $u^X_g: T_g(X)\to X$
is an isomorphism such that,
$$u^X_gT_g(u^X_h)= u^X_{gh}f^X_{g,h},\quad \mbox{for all\quad} g,h\in G.$$

(2)\quad A morphism $\phi: (Y,u_g^Y)\to (X,u_g^X)$ in $\mathcal{\C}^G$ is a morphism $\phi: Y\to X$ in $\mathcal{\C}$ such that $\phi u_g^Y=u_g^X\phi$, for all $g\in G$.

(3)\quad The tensor product in $\mathcal{\C}^G$ is  defined as $(Y,u_g^Y)\otimes (X,u_g^X)=(Y\otimes X, (u_g^Y\otimes u_g^X)j_g|_{Y,X})$, where $j_g|_{Y,X}:T_g(Y\otimes X)\to T_g(Y)\otimes T_g(X)$ is the isomorphism giving the monoidal structure on $T_g$.

\medbreak
There is a procedure opposite to equivariantization. Let $\C$ be  a fusion category and let $\Rep(G)\subseteq \mathcal{Z}(\C)$ be a Tannakian subcategory which embeds into $\C$ via the forgetful functor $\mathcal{Z}(\C)\to \C$, where $\mathcal{Z}(\C)$ is the Drinfeld center of $\C$. Let $A=\Fun(G)$ be the algebra of function on $G$. It is a commutative algebra in $\mathcal{Z}(\C)$ under the embedding  above. Let $\C_G$ be the category of left $A$-modules in $\C$. It is a fusion category which is called the de-equivariantization of $\C$ by $\Rep(G)$.

The two procedures above are inverse to each other; that is, there are canonical equivalences $(\C_G)^G\cong \C$  and $(\C^G)_G\cong \C$. Moreover, we have $$\FPdim(\C^G)=|G|\FPdim(\C) \mbox{\,\,and\,\,} \FPdim(\C_G)=\frac{\FPdim(\C)}{|G|}.$$

\medbreak
A braided $G$-crossed fusion
category is a fusion category $\C$ endowed with a $G$-grading $\C
= \oplus_{g \in G}\C_g$ and an action of $G$ by tensor autoequivalences
$\rho:\underline G \to \underline \Aut_{\otimes} \, \C$, such that $\rho^g(\C_h)
\subseteq
\C_{ghg^{-1}}$, for all $g, h \in G$, and a $G$-braiding $c: X \otimes Y \to
\rho^g(Y) \otimes X$, $g \in G$, $X \in \C_g$, $Y \in \C$, subject to appropriate
compatibility conditions.

\medbreak
Let $\C$ be a braided fusion category, and $\Rep(G)\subseteq \C$ be a Tannakian subcategory. The de-equivariantization $\C_G$ of $\C$ by $\Rep(G)$ is a braided $G$-crossed fusion category. The category $\C_G$ is not braided in general. But the neutral component $(\C_G)_e$ of the associated $G$-grading of $\C_G$ is braided. By \cite[Proposition 4.56]{drinfeld2010braided}, $\C$ is nondegenerate if and only if $(\C_G)_e$ is nondegenerate and the associated grading of $\C_G$ is faithful. In this case, $|G|^2$ divides $\FPdim(\C)$.

\section{Braided fusion categories of dimension $p^mq^nd$ and $p^2q^2r^2$}\label{sec3}
In this section, we study the braided fusion categories of dimension $p^mq^nd$ and $p^2q^2r^2$, where $p,q,r$ are distinct prime numbers, $d$ is a square-free natural number such that $(pq,d)=1$. Then we apply the results obtained to weakly integral braided fusion categories of dimension less than $1800$. Our results show all fusion categories involved are weakly group-theoretical.
\begin{proposition}\label{s_pqd}
Let $p,q$ be distinct prime numbers. Assume that $\C$ is a slightly degenerate fusion category of Frobenius-Perron dimension $p^mq^nd$, where $p,q$ are distinct prime numbers, $d$ is a square-free natural number such that $(pq,d)=1$. Then $\C$ is solvable.
\end{proposition}
\begin{proof}
By \cite[Corollary 3.4]{YU2020408}, $\frac{\FPdim(\C)}{2\FPdim(X)^2}$ is an algebraic integer for all $X\in \Irr(\C)$. Hence the Frobenius-Perron dimensions of integral simple objects have the form $p^aq^b$, where $a,b\geq 0$. Then $\C$ is weakly group-theoretical by \cite[Proposition 3.14]{YU2020408}. It follows from the arguments of \cite[Proposition 5.3]{Natale2018core} that every Tannakian subcategory of $\C$ is solvable if it exists. Furthermore, in this case $\C$ is solvable by \cite[Theorem 5.1]{Natale2018core}. In the rest of the paper, we aim to prove the existence of a nontrivial Tannakian subcategory.

We may assume that $\C$ contains no nontrivial non-degenerate fusion subcategories. In fact, if $\C$ contains such a fusion category $\D$ then $\C\cong \D\boxtimes \D'$ by Theorem \ref{MugerThm}, where $\D'$ is also slight degenerate. By induction, $\D'$ is solvable. By \cite[Corollary 5.4]{Natale2018core}, $\D$ is solvable. Hence $\C$ is solvable.

We assume on the contrary that $\C$ contains no nontrivial Tannakian subcategory. Then $G[X]=\{\1\}$ for all simple objects in $\C$ by \cite[Lemma 7.1]{natale2013weakly}. Moreover, every fusion subcategory of $\C$ is  slightly degenerate, also by \cite[Lemma 7.1]{natale2013weakly}.

Suppose first that $\C$ is integral. If the Frobenius-Perron dimensions of simple objects of $\C$ have a common prime factor $p$ or $q$ then the order of the group $G[X]$ is divisible by such a prime number. This is a contradiction. Therefore, there exists simple objects whose Frobenius-Perron dimension is a power of $p$ and also exists simple objects whose Frobenius-Perron dimension is a power of $q$. It follows from \cite[Proposition 7.4]{etingof2011weakly} that $\C$ contains a nontrivial Tannakian subcategory, a contradiction.

Suppose then that $\C$ is not integral. Then $\C$ is faithfully graded by an elementary abelian $2$-group whose trivial component is an integral fusion category $\C_e$, see \cite[Theorem 3.10]{gelaki2008nilpotent}. By the discussion above, $\C_e$ is also slightly degenerate. By the same arguments in the previous paragraph, we can prove that $\C_e$ contains a nontrivial Tannakian subcategory, also a contradiction.
\end{proof}

\begin{theorem}\label{mainth1}
Let $p,q$ be distinct prime numbers. Assume that $\C$ is a braided fusion category of Frobenius-Perron dimension $p^mq^nd$, where $d$ is a square-free natural number such that $(pq,d)=1$. Then $\C$ is weakly group-theoretical.
\end{theorem}
\begin{proof}
We may assume that $\C$ is neither non-degenerate nor slightly degenerate. Indeed, if $\C$ is non-degenerate or slightly degenerate, then $\C$ is solvable by \cite[Corollary 5.4]{Natale2018core} and Proposition \ref{s_pqd}.

Let $\E=\Rep(G)\subseteq \C'$ be the maximal Tannakian subcategory of $\C'$. Then the de-equivariantization $\C_G$ of $\C$ by $\Rep(G)$ is non-degenerate if $\E=\C'$, or slightly degenerate if $\E\varsubsetneq\C'$, see \cite[Corollary 4.31]{drinfeld2010braided}. Since the dimension of $\C_G$ still has the form $p^{m'}q^{n'}d'$, $\C_G$ is solvable by \cite[Corollary 5.4]{Natale2018core} and Proposition \ref{s_pqd}. Hence $\C$ is   weakly group-theoretical by \cite[Proposition 4.1]{etingof2011weakly}.
\end{proof}

\begin{remark}
In the proof above, we obtain that $\C_G$ is solvable, but we are not sure whether the group $G$ is solvable. Hence we can not determine the solvability of $\C$ . In fact, there exists a braided fusion category whose dimension has the form $p^mq^nd$, but it is not solvable. For example, let $\mathbb{A}_5$ be the alternating group of order $60=2^2\times 3\times 5$, and let $\Rep(\mathbb{A}_5)$ be the category of finite-dimensional representations of $\mathbb{A}_5$. Since $\mathbb{A}_5$ is a simple group, $\Rep(\mathbb{A}_5)$ is not solvable by \cite[Proposition 4.5]{etingof2011weakly}.
\end{remark}

\begin{corollary}\label{p2q2r2}
Let $p,q,r$ be distinct prime numbers. Assume that $\C$ is a braided fusion category of Frobenius-Perron dimension $p^2q^2r^2$. Then $\C$ is weakly group-theoretical.
\end{corollary}
\begin{proof}
Clearly, if $\C$ has a faithful grading then the trivial component matches the assumption of Theorem \ref{mainth1} and hence it is weakly group-theoretical. Thus $\C$ is weakly group-theoretical by \cite[Proposition 4.1]{etingof2011weakly}.

We may then assume that $\C=\C_{ad}$. In particular, $\C$ is integral. We shall prove that $\C$ contains a nontrivial Tannakian subcategory. In fact, if $\C$ contains a nontrivial Tannakian subcategory $\E=\Rep(G)$ then the de-equivariantizations $\C_G$ of $\C$ by $\Rep(G)$ is a braided $G$-crossed  fusion category whose trivial component $(\C_G)_e$ is braided and it is weakly group-theoretical by Theorem \ref{mainth1}. Therefore $\C$ is weakly group-theoretical by \cite[Proposition 4.1]{etingof2011weakly}.

If $\C$ is non-degenerate then the proof of \cite[Theorem 9.2]{etingof2011weakly} shows that $\C$ contains a nontrivial Tannakian subcategory.

If $\C$ is slightly degenerate then $\FPdim(\C)$ is even. We may assume that $p=2$ in this case. By \cite[Corollary 3.4]{YU2020408}, the dimensions of all simple objects of $\C$ are odd. Here we have two possibilities. The first possibility is that the dimensions of non-invertible simple objects $X$ are divisible by  $qr$. Then $qr$ divides the order of the group $G[X]$. By \cite[Lemma 2.4]{DONG2021386}, $\C_{ad}$ contains a nontrivial Tannakian subcategory. The second  possibility is that there exists simple objects whose Frobenius-Perron dimensions are power of $q$ or  $r$. It follows from \cite[Proposition 7.4]{etingof2011weakly} that $\C$ contains a nontrivial Tannakian subcategory.

If $\C$ is neither non-degenerate nor  slightly degenerate then the M\"uger  center $\C'$ contains a nontrivial Tannakian subcategory. This completes the proof.
\end{proof}

\begin{corollary}
Let $\C$ be a weakly integral braided fusion category of Frobenius-Perron dimension less than $1800$. Then $\C$ is weakly group-theoretical.
\end{corollary}

\begin{proof}
Let $n$ be a natural number less than $1800$. Then $n$ factorizes in the form $p^mq^nd$ or $p^2q^2r^2$, where $p,q,r$ are distinct prime numbers, $d$ is a square-free natural number such that $(pq,d)=1$. The result then follows from Theorem \ref{mainth1} and Corollary \ref{p2q2r2}.
\end{proof}

\begin{corollary}
Let $\C$ be a weakly integral braided fusion category. Assume that $\FPdim(\C)<33075$ is odd. Then $\C$ is solvable.
\end{corollary}

\begin{proof}
By Proposition \ref{wGt-sol}, it suffices to prove that $\C$ is weakly group-theoretical.

Let $n$ be an odd natural number less than $33075$. Then $n$ factorizes in the form $p^mq^nd$ except the case when $n=11025$, where $p,q$ are distinct prime numbers, $d$ is a square-free natural number such that $(pq,d)=1$. By Theorem \ref{mainth1}, it is enough to consider the case when $n=11025$. Notice that if $\FPdim(\C)=11025$ then $\C$ can not be slight degenerate. If $\C$ is non-degenerate then we are done by \cite[Theorem 8.2]{natale2013weakly}. If $\C$ is degenerate then the  M\"uger  center $\C'$ is a nontrivial Tannakian subcategory $\E=\Rep(G)$. Then the de-equivariantizations $\C_G$ of $\C$ by $\Rep(G)$ is a braided $G$-crossed  fusion category whose trivial component $\C_G^0$ is braided and it is weakly group-theoretical. Hence $\C$  is weakly group-theoretical.
\end{proof}

\section{Non-braided fusion categories of dimension $84$ and $90$}\label{sec4}
We start this section with the existence of nontrivial symmetric categories in the Drinfeld center $\Z(\C)$.

\subsection{Existence of nontrivial symmetric categories}\label{sec41}
There is an obvious forgetful tensor functor $F:\Z(\C)\to \C$. By \cite[Proposition 3.39]{etingof2004finite}, the forgetful functor $F:\Z(\C)\to \C$ is surjective.

\begin{lemma}\cite[Lemma 2.1]{dong2012frobenius}\label{lemma1}
Let $F_0 : G(\Z(\C)) \to G(\C)$ be the group homomorphism induced by the forgetful
functor $F : \Z(\C)\to \C$. Then the following hold:

(1)\, $\C$ is faithfully graded by the group $\widehat{N}$, where $N$ is the kernel of $F_0$.

(2)\, Suppose $\U(\C)$ is trivial. Then the group homomorphism $F_0$ is injective.
\end{lemma}

Since $\Z(\C)$ is semisimple the forgetful functor $F:\Z(\C)\to \C$ has a right adjoint $I: \C\to \Z(\C)$. By \cite[Lemma 3.2]{etingof2011weakly}, $I(\1)\in \Z(\C)$ has a natural structure of commutative algebra. By \cite[Proposition 5.4]{etingof2005fusion}, $\FPdim(I(\1))=\FPdim(\C)$.

\begin{lemma}\label{lemma2}
Assume that $I(\1)$ contains nontrivial invertible simple objects of $\Z(\C)$. Then $\C$ is faithful graded by some finite group.
\end{lemma}
\begin{proof}
Let $g\in I(\1)$ be a nontrivial invertible object. Then $\Hom_{\C}(g,I(\1))=\Hom_{\C}(F(g),\1)\neq 0$ implies that $F(g)=\1$. Hence the kernel of $F_0 : G(\Z(\C)) \to G(\C)$ is not trivial. By Lemma \ref{lemma1}, $\C$ is faithfully graded by the group $\widehat{N}$, where $N$ is the kernel of $F_0$.
\end{proof}

 The following lemma is contained in the proof of \cite[Lemma 9.17]{etingof2011weakly}. We explicitly state it for reader's convenience.

\begin{lemma}\label{lemma3}
Let $\D\subset \C$ be a fusion subcategory. Then $I(\1)$ contain a subalgebra $B$ corresponding to $\D$ such that $\FPdim(B)=\FPdim(\C)/\FPdim(\D)$.
\end{lemma}

\begin{lemma}\label{lemma4}
Assume that the Drinfeld center $\Z(\C)$ contains a nontrivial symmetric category $\E$. If the order of $G(\C)$ is odd then $\Z(\C)$ contains a nontrivial Tannakian subcategory.
\end{lemma}
\begin{proof}
We may assume that $\C$ has a trivial universal grading, otherwise \cite[Proposition 2.9(ii)]{etingof2011weakly} shows that $\Z(\C)$ contains a nontrivial Tannakian subcategory.

Consider the group homomorphism $F_0:G(\Z(\C)))\to G(\C)$. Since $\U(\C)$ is trivial, $F_0$ is injective by Lemma \ref{lemma1}. It follows that the order of $G(\Z(\C))$ is odd since the order of  $G(\C)$ is odd. Hence,  $\Z(\C)$ can not contain fusion subcategory of dimension $2$. It follows that the dimension of  $\E$ is greater that $2$. By Lemma \ref{lem22}, $\E$ contains a nontrivial Tannakian subcategory.
\end{proof}

Let $\C$ be a pre-modular category with a spherical structure $\psi$. Let $\Tr$ denote the trace corresponding to $\psi$. The S-matrix $S = \{s_{X,Y} \}_{X,Y\in \Irr(\C)}$ of $\C$ is defined by $s_{X,Y}=\Tr(c_{Y,X}c_{X,Y})$.

\begin{theorem}\label{dim}
Let $p,q,r$ be distinct prime numbers and $n$ be a positive integer. Assume that $\C$ is an integral non-degenerate fusion category of dimension $p^{2n}q^2r^2$. Then one of the following holds true.

(1)\, $\C$ contains a nontrivial symmetric subcategory.

(2)\, The Frobenius-Perron dimensions of simple objects of $\C$ can not be divisible by $p^n$.
\end{theorem}
\begin{proof}
Let $X$ be a simple object of $\C$. Then $\FPdim(X)$ divides $p^nqr$ by \cite[Theorem 2.11]{etingof2011weakly}.  If $\FPdim(X)$ is a power of $p,q$ or $r$ then $\C$  contains a nontrivial symmetric category by \cite[Corollary 7.2]{etingof2011weakly}. This proves (1).

In the rest of the proof, we assume that $\FPdim(X)$ is not a power of a prime number and prove $\FPdim(X)$ can not be $p^nd$, where $d=1,q$ or $r$ .

We first consider the case when $\C$ does not contain nontrivial invertible simple objects. In this case, $\C$ must contain a simple object of dimension $qr$. In fact, if not, all possible dimensions of nontrivial simple objects have a common factor $p$ which implies that $\C$ has a pointed fusion subcategory with dimension divisible  by $p$. This contradicts our assumption.

Assume on the contrary there exists a simple object $X_0$  of dimension $p^nd$. By the orthogonality of columns of the S-matrix, we have:

\begin{equation}\label{eq1}
\begin{split}
\sum_{X\in \Irr{(\C)}}\frac{s_{X_0,X}}{\FPdim(X_0)}\FPdim(X)=0.
\end{split}
\end{equation}
Hence there exists $X_1\in \Irr(\C)$ of dimension $qr$ such that $s_{X_0,X_1}\neq 0$. In fact, if there does not exist such a simple object, then the simple object $X$ such that $s_{X_0,X}\neq 0$ is either $1$ or its dimension has a prime factor $p$. It follows that the left side of the equation (\ref{eq1}) is equal to $1$ modular $p$, a contradiction.

The ratios $\frac{s_{X_0,X_1}}{\FPdim(X_0)}$ and $\frac{s_{X_0,X_1}}{\FPdim(X_1)}$ are both algebraic integers, and so is $\frac{s_{X_0,X_1}}{p^nqr}$ . This implies that $\frac{s_{X_0,X_1}}{\FPdim(X_0)}$ is divisible by $t:=\frac{qr}{d}$.

On the other hand, we have
\begin{equation}\label{eq3}
\begin{split}
\sum_{X\in \Irr{\C)}}\left |\frac{s_{X_0,X}}{\FPdim(X_0)}\right |^2=\frac{\FPdim(\C)}{\FPdim(X_0)^2}=t^2.
\end{split}
\end{equation}

As we know, $s_{X_0,X}$ is a sum of roots of unity, hence every summand on the left side is a totally positive algebraic integer. The summand corresponding to $X=\1$ is $1$ and the summand $a$ corresponding to $X=X_1$ is an algebraic integer disvisble by $t^2$. Hence there is a Galois automorphism $\delta$ such that $\delta(a)\geq t^2$. Applying $\delta$ to equation (\ref{eq3}), we get that the left side $\geq 1+t^2$, a contradiction.

Now we consider the case when $\C$ has nontrivial invertible simple objects. Let $\B=\C_{pt}$ be the maximal pointed fusion subcategory of $\C$. If $\B$ is degenerate then the M\"uger center of $\B$ is a nontrivial symmetric subcategory.

If $\B$ is non-degenerate then $\C=\B\boxtimes \B'$ by Theorem \ref{MugerThm}, where $\B'$ is the M\"uger centeralizer of $\B$ in $\C$. In particular,  $\B'$ is non-degenerate and does not contain nontrivial invertible simple objects. By the proof of the first case, the Frobenius-Perron dimensions of simple objects of $\B'$ and hence $\C$ can not be $p^nd$, where $d=1,q$ or $r$. This completes the proof.
\end{proof}

\begin{proposition}\label{eXt_Tann}
Let $\C$ be a fusion category of dimension $2qp^n$, where $p,q$ are distinct odd primes and $n\geq 1$. If $\Z(\C)$ contains a nontrivial symmetric subcategory then $\Z(\C)$ contains a nontrivial Tannakian subcategory.
\end{proposition}
\begin{proof}
Assume on the contrary that  $\Z(\C)$ has a unique nontrivial symmetric subcategory $\E$ which is equivalent to the category of super vectors. Let $\A=\E'\subset \Z(\C)$ be the M\"uger centralizer of $\E$ in $\Z(\C)$. Then the M\"uger center $\Z_2(\A)$of $\A$ is $\E$ and hence $\A$ is a slightly degenerate fusion category of dimension $2 q^2 p^{2n}$.

By \cite[Proposition 7.4]{etingof2011weakly}, we may assume that $\A$ can not contain simple objects of odd prime power dimension. On the other hand,  \cite[Corollary 3.4]{YU2020408} shows that $\frac{\FPdim(\A)}{2\FPdim(X)^2}$ is an integer for all noninvertible simple objects $X$ of $\A$. Hence $\FPdim(X)=p^iq^j$ for $i,j\geq 1$. The decomposition of $X\otimes X^*\in\Z(\C)_{ad}$ implies that the order of $G[X]$ is divisible by $pq$. By \cite[Lemma 2.4]{DONG2021386}, the fusion subcategory generated by  $G[X]$ is a symmetric subcategory which must be Tannakian because it has odd dimension. This is a contradiction.
\end{proof}

Let $X,Y$ be simple objects of a braided fusion category. We say that $X,Y$ projectively centralize each other if $c_{Y,X}c_{X,Y} =\lambda \id_{X\otimes Y}$ for some $\lambda\in k^{\times}$.

\begin{proposition}\label{slig-deg-exist-Tanna}
Let $\C$ be a slightly degenerate fusion category. Suppose that $\C$ contains a simple object of dimension $2$ or $4$. Then $\C$ contains a nontrivial Tannakian subcategory.
\end{proposition}
\begin{proof}
Let $\A$ be the category generated by the invertible objects of $\C$. By \cite[Proposition 2.6(ii)]{etingof2011weakly}, $\A=\svect\boxtimes \A_0$, where $\svect$ is the M\"uger center of $\C$ and $\A_0$ is a non-degenerate pointed category. Then $\C=\A_0\boxtimes \A_0'$ by Theorem \ref{MugerThm}, where $\A_0'$ is the M\"uger centralizer of $\A_0$. In particular, $\A_0'$ is slightly degenerate and contains only two invertible objects: $\1$ and the generator $\delta$ of $\svect$. From $\C=\A_0\boxtimes \A_0'$, we know $\A_0'$  also contains a simple object $X$ of dimension $2$ or $4$. By \cite[Proposition 2.6(i)]{etingof2011weakly}, $\delta$ can not appear in the decomposition of $X\otimes X^*$.

 If $\FPdim(X)=2$ then we must have $X\otimes X^*=\1\oplus Y$, where $Y$ is a simple object of dimension $3$. By \cite[Proposition 7.4]{etingof2011weakly},  $\A_0'$  contains a nontrivial Tannakian subcategory.

 If $\FPdim(X)=4$ then we must have $X\otimes X^*=\1\oplus Y$, where $Y$ is a direct sum of noninvertible simple objects. Since $\FPdim(Y)=15$, $Y$ at least contains a simple object of odd dimension. If there is an object of dimension $3,5,7,9,11$ or $13$, then we are done by \cite[Proposition 7.4]{etingof2011weakly}. So it suffices to consider the case when $15$ is the only possible dimension. Thus we have $X\otimes X^*=\1\oplus Y$, where $Y$ is a simple object of dimension $15$. Then we have $s_{X,X^*}=\alpha+15\beta$, where $\alpha,\beta$ are roots of unity. Meanwhile, $s_{X,X^*}$ is divisible by $\FPdim(X)=4$. Hence $\alpha-\beta$ is divisible by $4$ which implies that $\alpha-\beta=0$ and hence $s_{X,X^*}=16\alpha$. Thus $X$ projectively centralizes its dual and hence itself. So $Y$ centralizes itself. Then the subcategory generated by $Y$ is a symmetric category which must contain a nontrivial Tannakian subcategory.
\end{proof}

\begin{proposition}\label{lemma5}
Let $p,q,r$ be distinct prime numbers, and let $\C$ be a  fusion category of dimension $p^2qr$. Assume that the Drinfeld  center $\Z(\C)$ contains a nontrivial Tannakian subcategory $\E=\Rep(G)$. Then

(1) If $\FPdim (\E)<p^2qr$ then $\C$ is solvable.

(2) If $\FPdim (\E)=p^2qr$ then $\C$ is group-theoretical.
\end{proposition}

\begin{proof}
Assume that  $\FPdim (\E)$ is a power of $p$. Then $G$ is a solvable group. It follows that $G$ has a quotient group $H$ such that $|H|$ is $p$. Hence $\E$ has a subcategory $\Rep(H)$. Under the forgetful functor $F:\Z(\C)\to \C$, $\Rep(H)$ either maps to $\vect$, or embeds to $\C$. Hence $\C$ is an $H$-extension of some fusion category $\D_1$  by \cite[Proposition 2.9(i)]{etingof2011weakly}, or $\C$ is an $H$-equivariantization of some fusion category $\D_2$  by \cite[Proposition 2.10(i)]{etingof2011weakly}. In both cases, $\FPdim (\D_1)=\FPdim(\D_2)=pqr$. By Corollary \ref{dimpqr}, $\D_1$ and $D_2$ are both solvable, hence $\C$ is solvable by \cite[Proposition 4.5]{etingof2011weakly}.

\medbreak
Assume that $\FPdim (\E)$ has a prime factor $q$ or $r$. We consider the de-equivariantization $\Z(\C)_G$ of $\Z(\C)$ by $\E$. Set $\D=\Z(\C)_G$. Then $\D=\oplus_{g\in G}\D_g$ is faithfully graded by $G$, see \cite[Proposition 4.56]{drinfeld2010braided}. Hence $\FPdim (\D_e)=\frac{p^{4}q^2r^2}{\FPdim(\E)^2}$.

If $\FPdim(\E)<p^2qr$ then either $\FPdim(\E)=pqr$  or $\FPdim(\E)$ has at most $2$ distinct prime factors. This implies that $G$ is a solvable group. In this case, $\FPdim (\D_e)=pqr$ or $\FPdim (\D_e)$ has at most $2$ distinct prime factors. By Corollary \ref{dimpqr} and \cite[Theorem 1.6]{etingof2011weakly}, $\D_e$ is solvable in both cases. Therefore, $\Z(\C)$ and hence $\C$ are solvable by \cite[Proposition 4.5]{etingof2011weakly}.

If $\FPdim(\E)=p^2qr$ then $\FPdim (\D_e)=1$ and hence $\D$ is pointed. It follows that $\Z(\C)$, being an equivariantization of a pointed fusion category, is group-theoretical by \cite[Theorem 7.2]{naidu2009fusion}. So $\C$ is group-theoretical by \cite[Theorem 1.5]{drinfeld2007g}.
\end{proof}

\begin{corollary}\label{cordimp2qr}
Let $p,q,r$ be distinct prime numbers. Assume that $\C$ is a  weakly group-theoretical fusion category of dimension $p^2qr$. Then $\C$ is group-theoretical or solvable.
\end{corollary}
\begin{proof}
By \cite[Proposition 4.2]{etingof2011weakly}, $\Z(\C)$ contains a nontrivial Tannakian subcategory. The result then follows from Proposition \ref{lemma5}.
\end{proof}

\subsection{Fusion categories of dimension $84$}\label{sec42}

For any fusion category $\C$, we use $\deg(\C)$ to denote the set of Frobenius-Perron dimensions of simple objects of $\C$.

\begin{lemma}\label{exist_inv_simple}
Let $\C$ be a fusion category. Assume $\deg(\C)=\{1,6,14,21,42\}$ and $G[X]$ is trivial for all $X\in\Irr(\C)$. Then $\C$ has at least $6$ nontrivial invertible simple objects. Moreover, if $\C$ is the Drinfeld center of a fusion category $\D$ and the universal grading group $\U(\D)$ is trivial, then $\D$ also has at least $6$ nontrivial invertible simple objects.
\end{lemma}
\begin{proof}
Let $X_6$ be a simple object of dimension $6$. Then $X_6\otimes X_6^*=\1+X_{14}+X_{21}$, where $X_{14},X_{21}$ are simple objects of dimension $14$ and $21$, respectively. From $m(X_{14},X_6\otimes X_6^*)=m(X_6,X_{14}\otimes X_6)=1$, we can write $X_{14}\otimes X_6=X_6+W$, where $W$ is a direct sum of simple objects of dimension $6,14,21$ or $42$ and $\FPdim (W)=78$. Then we have an equation $6a_1+14a_2+21a_3+42a_4=78$ which show that $a_1=6$ or $8$. This implies that there are $6$ or $8$ $6$-dimensional simple objects in the decomposition of $W$.

Let $X_6'$ be a simple object of dimension $6$ such that $m(X_6',X_{14}\otimes X_6)\geq1$. Then $m(X_6',X_{14}\otimes X_6)=m(X_{14},X_6'\otimes X_6^*)\leq 2$. If $m(X_{14},X_6'\otimes X_6^*)=1$ then $X_6'\otimes X_6^*=X_{14}+W$ where $\FPdim (W)=22$. Considering the possible decomposition of $W$,  we get that $W$ at least contains one invertible simple object.
If $m(X_{14},X_6'\otimes X_6^*)=2$ then $X_6'\otimes X_6^*=2X_{14}+W$ where $\FPdim (W)=8$. Considering the possible decomposition of $W$,  we get that $W$ at least contains two invertible simple objects. It is easy to check all invertible objects obtained above are pairwise different. It follows that $\C$ has at least $6$ nontrivial invertible simple objects.

\medbreak

Assume that $\C=\Z(\D)$. Since $\U(\D)$ is trivial, the group homorphism $F_0:G(\C)\to G(\D)$ is injective by Lemma \ref{lemma1}. Hence $\D$ also has at least $6$ nontrivial invertible simple objects.
\end{proof}

\begin{lemma}\label{exist_symm1}
Let $\C$ be an integral fusion category of dimension $84$. Assume that $\C$ has a subcategory $\B$ of dimension $\geq 7$. Then the Drinfeld center $\Z(\C)$ of $\C$ contains a nontrivial symmetric subcategory.
\end{lemma}
\begin{proof}
By Lemma \ref{lemma3}, there exists a subalgebra $D\subset I(\1)$ corresponding to the fusion subalgebra $\B$ such that $\FPdim(D)=\FPdim(\C)/\FPdim(\B)\leq 12$. In view of \cite[Theorem 2.11]{etingof2011weakly}, the Frobenius-Perron dimensions of simple objects of $\Z(\C)$ divides $84$. The possible decompositions of $D$ as an object of $\Z(\C)$ shows that $D$ either contains nontrivial invertible simple objects, or contains simple objects of prime power dimension. If the former case holds true then $\C$ has a nontrivial grading by Lemma \ref{lemma2} and hence $\Z(\C)$ contains a nontrivial Tannakian subcategory by \cite[Corollary 7.2]{etingof2011weakly}. If the later case holds true then $\Z(\C)$  contains a nontrivial symmetric subcategory by \cite[Corollary 7.2]{etingof2011weakly}.
\end{proof}

\begin{lemma}\label{exist_symm2}
Let $\C$ be an integral fusion category of dimension $84$. Then the Drinfeld center $\Z(\C)$ of $\C$ contains a nontrivial symmetric subcategory.
\end{lemma}
\begin{proof}
By Theorem \ref{dim}, we may assume that the Frobenius-Perron dimensions of simple objects of $\Z(\C)$ are $1,6,14,21$ or $42$. If there exists $X\in\Irr(\Z(\C))$ such that $G[X]$ is not trivial then $G[X]\subset\Z(\C)_{ad}\cap\Z(\C)_{pt}$ and hence it is a symmetric subcategory by \cite[Lemma 2.4]{DONG2021386}. We are done in this case. This fact also implies that  we may assume $\Z(\C)$ has  simple objects of dimension $6$, otherwise the order of $G[X]$ is divisible by $7$ for any noninvertible simple object $X\in\Irr(\Z(\C))$. This is because the Frobenius-Perron dimensions of noninvertible simple objects have a common prime factor $7$.

Now we can assume that $\deg(\Z(\C))=\{ 1,6,14,21, 42\}$ and $G[X]$ is trivial for any noninvertible simple object $X$. In addition, we also assume that $\U(\C)$ is trivial, otherwise \cite[Proposition 2.9(ii)]{etingof2011weakly} shows that $\Z(\C)$ contains a nontrivial Tannakian subcategory. Therefore, we match the assumption of Lemma \ref{exist_inv_simple}. Hence $G(\C)$ has at least $7$ invertible simple objects. It follows from Lemma \ref{exist_symm1} that  $\Z(\C)$ contains a nontrivial symmetric subcategory.
\end{proof}

\begin{theorem}\label{main1}
Let $\C$ be a fusion category of dimension $84$. Then $\C$ is solvable or group-theoretical.
\end{theorem}
\begin{proof}
If $\C$ is not integral then $\C$ is a $G$-extension of an integral fusion subcategory $\D$ with lower dimension by \cite[Theorem 3.10]{gelaki2008nilpotent}, where $G$ is an elementary abelian $2$-group. Then $\FPdim(\D)=21$ or $42$. By  Corollary \ref{dimpqr} and \cite[Theorems 1.6]{etingof2011weakly}, $\D$ is solvable. Hence $\C$ is solvable by \cite[Proposition 4.5]{etingof2011weakly}.

We therefore assume that $\C$ is integral. By Lemma \ref{exist_symm2},  $\Z(\C)$ has a nontrivial symmetric subcategory. In view of Lemma \ref{lem22}, we may assume that $\Z(\C)$ has a unique nontrivial symmetric subcategory $\E$ which is equivalent to the category of super vectors. Let $\B=\E'\subset \Z(\C)$ be the M\"uger centralizer of $\E$ in $\Z(\C)$. Then the M\"uger center $\Z_2(\B)$of $\B$ is $\E$ and hence $\B$ is a slightly degenerate fusion category.

Let $\A$ be the category spanned by the invertible objects of $\B$. By \cite[Proposition 2.6(ii)]{etingof2011weakly}, $\A=\E\boxtimes \A_0$, where $\A_0$ is a non-degenerate pointed category. Then $\B=A_0\boxtimes \A_0'$ by Theorem \ref{MugerThm}. It is known that $\A_0'$ is slightly degenerate fusion category and contains only two invertible objects: $\1$ and the generator $\delta$ of $\E$.

We claim that $\A_0'$ contains a simple object $X$ of prime power dimension. In fact, if $\A_0'$ does not contain simple objects of prime dimension then $\deg(\A_0')=\{1,16,14,21,42\}$. Then  $\A_0'$ contains at least $7$ invertible simple objects by Lemma \ref{exist_inv_simple} which contradicts the fact $G(\A_0')=\{\1,\delta\}$. If $\FPdim(X)$ is odd then $\A_0'$ has a nontrivial Tannakian subcategory by \cite[Proposition 7.4]{etingof2011weakly}. If $\FPdim(X)$ is even then $\A_0'$ has a nontrivial Tannakian subcategory by Proposition \ref{slig-deg-exist-Tanna}.  In both case, $\Z(\C)$ has a nontrivial Tannakian subcategory. The result then follows from Lemma \ref{lemma5}.
\end{proof}

\subsection{Fusion categories of dimension $90$}\label{sec43}
\begin{lemma}\label{lemma11}
Let $\C$ be an integral fusion category of dimension $90$. If $\C$ has a fusion subcategory $\D$ of dimension $\geq 6$. Then one of following holds true.

(1) $\C$ is faithfully graded by a nontrivial finite group.

(2) $\Z(\C)$ contains a nontrivial symmetric subcategory. In particular, $\C$ is weakly group-theoretical.
\end{lemma}
\begin{proof}
Let $B\subseteq A=I(\mathbf{1})$ be the subalgebra of $A$ corresponding to the fusion subcategory $\D$. Then $\FPdim(B)=\frac{90}{\FPdim(\D)}\leq 15$. In view of [6,Theorem 2.11], the  Frobenius-Perron dimension of simple object of $\Z(\C)$ divide 90. Hence the possible decomposition of $B$ as an object of $\Z(\C)$ shows that $B$ (and hence $A$) contains nontrivial invertible simple object or simple object of prime power dimension. If the former possibility holds then $\C$ is faithfully graded by a nontrivial finite group by Lemma \ref{lemma2}. In this case, the trivial component has dimension $\leq 45$ and hence is weakly  group-theoretical. Thus $\C$ is weakly group-theoretical.

If the later possibility holds then $Z(\C)$ contains a nontrivial symmetric subcategory by \cite[Corollary 7.2]{etingof2011weakly} . In this case, $\C$ is weakly group-theoretical by Proposition \ref{eXt_Tann} and Lemma \ref{lemma5}.
\end{proof}

\begin{theorem}\label{theorem90}
Let $\C$ be a fusion category of dimension 90. Then $\C$ is weakly group-theoretical.
\end{theorem}
\begin{proof}
If $\C$ is not integral then $\C$ is a $\mathbb{Z}_2$-extension of an integral fusion subcategory $\D$ with dimension $45$, by \cite[Theorem 3.10]{gelaki2008nilpotent}. By \cite[Theorem 1.6]{etingof2011weakly}, $\D$ is solvable. Hence $\C$ is weakly group-theoretical by \cite[Proposition 4.1]{etingof2011weakly}.

\medbreak
If $\C$ is integral then it suffices to prove that $\C$ has a fusion subcategory of dimension $\geq 6$, by Lemma \ref{lemma11}.

By [Theorem 4.2], $\C$ has the following possible types:
\begin{align*}& (1, 2; 2, 4; 3, 8), (1, 2; 2, 4; 6, 2), (1, 2; 2, 4; 3, 4;  6, 1), (1, 2; 2,
22),&\\
& (1, 6; 2, 3; 3, 8), (1, 6; 2, 3; 6, 2), (1, 6; 2, 3; 3, 4; 6, 1), (1, 6; 2, 21), &\\
& (1, 9; 3, 9), (1, 9; 3, 1; 6, 2), (1, 9; 3, 5; 6, 1), (1, 9; 9, 1), (1, 10; 2, 20), (1, 15; 5, 3),& \\
&
 (1, 18; 2, 18), (1, 18; 3, 8), (1, 18; 3, 4; 6, 1), (1, 18; 6, 2),  (1, 30; 2, 15), (1, 45; 3, 5).&
\end{align*}

Assume that $\C$ has the first $3$ types. Then $G[X]=G(\C)$ for any simple object $X$ of dimension $2$. In fact, if there exists $2$-dimensional simple object $X$ such that $G[X]=\{1\}$ then Theorem \ref{thm-dim-2} shows that $\C$ has a fusion subcategory of dimension $12, 24$ or $60$. It is  impossible since they can not divides $90$. It follows from \cite[Lemma 3.2]{dong2012frobenius} that all simple object of dimension $1$ and $2$ generate a fusion subcategory $\D$ of dimension $18$.

Assume that $\C$ has the forth type. Then $X\otimes X^*=1+g+X_2$ for every 2-dimensional simple object $X$, where $G(\C)=\{1,g\}$ and $X_2$ is $2$-dimensional simple object. Hence $X_2$ is a self-dual simple object. Let $\C(X_2)$ be the fusion subcategory generated by $X_2$. If  $\C\langle X_2\rangle$ is a proper subcategory of $\C$ then $\FPdim(\C(X_2))\geq6$. If $\C\langle X_2\rangle$ then $\C$ is Grothendieck equivalent to the category $\Rep(\D_{45})$ by \cite[Theorem 1.2]{WangDongLi2022}, where $\D_{45}$ is the dihedral group of order $90$. Hence $\C$ is group-theoretical by \cite[Proposition 4.7(1)]{naidu2011finiteness}.

For the remaining types, the largest pointed fusion subcategory $\C_{pt}$ has dimension $\geq 6$. This completes the proof.
\end{proof}

Combing Proposition \ref{wGt-sol} and Theorem \ref{theorem90}, we get the following corollary.
\begin{corollary}\label{sol-90}
Let $\C$ be a fusion category of dimension $90$. Then $\C$ is solvable.
\end{corollary}

\subsection{Main Results}\label{sec44}
\begin{theorem}\label{dim120}
Let $\C$ be a weakly integral fusion category of dimension less than $120$. Then $\C$ is either group-theoretical or solvable.
\end{theorem}
\begin{proof} Suppose $\FPdim(\C) = n$ is a natural number
and  $n< 120$.

If $n=p^aq^b$, where $p$, $q$ are prime numbers, $a, b \geq 0$, then $\C$ is
solvable by \cite[Theorem 1.6]{etingof2011weakly}.

If $n=84$ or $90$ then $\C$ is solvable by Theorem \ref {main1} and Corollary \ref{sol-90}.

If $n = pqr$, where $p$, $q$ and $r$ are distinct prime numbers, then  $\C$ is solvable by Corollary \ref{dimpqr}.

It remains to consider the case when $n=60$. We shall follow the line of proof of Proposition \ref{wGt-sol}.
By \cite[Theorems 9.16]{etingof2011weakly}, $\C$ is weakly group-theoretical and hence it is Morita equivalent to a nilpotent fusion category $\D$. If $\D$ is pointed then $\C$ is group-theoretical.
\medbreak
We then assume that $\D$ is not pointed and hence if there is a sequence of fusion subcategories
\begin{equation}\label{nilp1}
\begin{split}
\mathcal{\D}_0 ={\rm Vec}, \mathcal{\D}_1, \cdots, \mathcal{\D}_n =\mathcal{\D}
\end{split}
\end{equation}
and a sequence of finite groups $G_1,\cdots, G_n$ such that $\mathcal{\D}_i$ is obtained from $\mathcal{\D}_{i-1}$ by a $G_i$-extension.

For every group $G_i$, $1\leq i\leq n$, the order $|G_i|$ divides $60$ and hence it is solvable. Moreover, the fusion subcategory $\D_1$ must be pointed and hence $\D_1=\vect_{K,\omega}$,  for some $3$-cocycle $\omega\in H^3(G,k^{\times})$ and a solvable group $K$.  Then $\C$ is solvable by the same reason as in the proof of Proposition \ref{wGt-sol}.
\end{proof}

\begin{corollary}
Let $\C$ be a strictly weakly integral fusion category of dimension less than $120$. Then $\C$ is solvable.
\end{corollary}
\begin{proof}
The result follows from Theorem \ref{dim120} and the fact that a group-theoretical fusion category is integral, see \cite[Corollary 8.43]{etingof2005fusion}.
\end{proof}

By Theorem \ref{dim120} and \cite[Theorems 1.5]{etingof2011weakly}, we can recover the main result in \cite{dong2012frobenius}.
\begin{corollary}
Let $\C$ be a weakly integral fusion category of dimension less than $120$. Then $\C$  has the strong Frobenius property.
\end{corollary}

\section*{Acknowledgements}
The research of the author is supported by the Natural Science Foundation of
Jiangsu Providence (Grant No. BK20201390).



\end{document}